\documentclass[12pt, hidelinks]{article}
\usepackage[small, bf]{caption}
\usepackage[utf8]{inputenc}
\usepackage{fullpage}
\usepackage{appendix}
\usepackage{graphicx}
\usepackage{amsmath}
\usepackage{amssymb}
\usepackage{amsthm, amsfonts}
\usepackage{url}
\usepackage{subfig}

\usepackage{xcolor}
\usepackage{listings}

\definecolor{darkred}{rgb}{0.6,0.0,0.0}
\definecolor{darkgreen}{rgb}{0,0.50,0}
\definecolor{lightblue}{rgb}{0.0,0.42,0.91}
\definecolor{orange}{rgb}{0.99,0.48,0.13}
\definecolor{grass}{rgb}{0.18,0.80,0.18}
\definecolor{pink}{rgb}{0.97,0.15,0.45}

\lstdefinestyle{colored}{ %
  basicstyle=\ttfamily,
  backgroundcolor=\color{white},
  commentstyle=\color{green}\itshape,
  keywordstyle=\color{blue}\bfseries\itshape,
  stringstyle=\color{red},
}

\lstdefinelanguage{PythonPlus}[]{Python}{
  morekeywords=[1]{,as,assert,nonlocal,with,yield,self,True,False,None,} 
  morekeywords=[2]{,__init__,__add__,__mul__,__div__,__sub__,__call__,__getitem__,__setitem__,__eq__,__ne__,__nonzero__,__rmul__,__radd__,__repr__,__str__,__get__,__truediv__,__pow__,__name__,__future__,__all__,}, 
  morekeywords=[3]{,object,type,isinstance,copy,deepcopy,zip,enumerate,reversed,list,set,len,dict,tuple,range,xrange,append,execfile,real,imag,reduce,str,repr,}, 
  morekeywords=[4]{,Exception,NameError,IndexError,SyntaxError,TypeError,ValueError,OverflowError,ZeroDivisionError,}, 
  morekeywords=[5]{,ode,fsolve,sqrt,exp,sin,cos,arctan,arctan2,arccos,pi, array,norm,solve,dot,arange,isscalar,max,sum,flatten,shape,reshape,find,any,all,abs,plot,linspace,legend,quad,polyval,polyfit,hstack,concatenate,vstack,column_stack,empty,zeros,ones,rand,vander,grid,pcolor,eig,eigs,eigvals,svd,qr,tan,det,logspace,roll,min,mean,cumsum,cumprod,diff,vectorize,lstsq,cla,eye,xlabel,ylabel,squeeze,}, 
}

\lstdefinestyle{colorEX}{
  basicstyle=\footnotesize\ttfamily,
  backgroundcolor=\color{white},
  commentstyle=\color{darkgreen}\slshape,
  keywordstyle=\color{blue}\bfseries\itshape,
  keywordstyle=[2]\color{blue}\bfseries,
  keywordstyle=[3]\color{grass},
  keywordstyle=[4]\color{red},
  keywordstyle=[5]\color{orange},
  stringstyle=\color{darkred},
  emphstyle=\color{pink}\underbar,
}

\usepackage{tikz}
\usetikzlibrary{positioning,chains,fit,shapes,calc}
\usepackage{listings}

\let\phi\varphi

\newcommand{\dmax}{\Delta^\mathrm{max}}
\newcommand{\ones}{\mathbf 1}
\newcommand{\reals}{{\mbox{\bf R}}}

\newcommand{\eg}{{\it e.g.}}
\newcommand{\ie}{{\it i.e.}}

\newcommand{\BEAS}{\begin{eqnarray*}}
\newcommand{\EEAS}{\end{eqnarray*}}
\newcommand{\BEA}{\begin{eqnarray}}
\newcommand{\EEA}{\end{eqnarray}}
\newcommand{\BEQ}{\begin{equation}}
\newcommand{\EEQ}{\end{equation}}
\newcommand{\BIT}{\begin{itemize}}
\newcommand{\EIT}{\end{itemize}}

\title{Optimal Routing for Constant Function Market Makers}

\author{
Guillermo Angeris\\
{\small \texttt{angeris@stanford.edu}}
\and
Tarun Chitra\\
{\small \texttt{tarun@gauntlet.network}}
\and
Alex Evans\\
{\small \texttt{aevans@baincapital.com}}
\and
Stephen Boyd\\
{\small \texttt{boyd@stanford.edu}}}

\date{December 2021}

\begin{document}

\maketitle

\begin{abstract}
We consider the problem of optimally executing an order involving multiple
crypto-assets, sometimes called tokens,
on a network of multiple constant function market 
makers (CFMMs).
When we ignore the fixed cost associated with executing an order on a 
CFMM,
this optimal routing problem can be cast as a convex optimization problem, 
which is computationally tractable.
When we include the fixed costs, the optimal routing problem
is a mixed-integer convex problem, which can be solved
using (sometimes slow) global optimization methods, or 
approximately solved using various heuristics based on 
convex optimization.
The optimal routing problem includes as a special case the problem 
of identifying an arbitrage present in a network of CFMMs, or certifying
that none exists.
\end{abstract}

\section*{Introduction}
Decentralized exchanges (DEXs) are a popular application of public blockchains that allow users to
trade assets without the need for a trusted intermediary to facilitate the exchange. DEXs are
typically implemented as \emph{constant function market makers}
(CFMMs)~\cite{angerisImprovedPriceOracles2020}. In CFMMs, liquidity providers contribute reserves of
assets. Users can then trade against these reserves by tendering baskets of assets in exchange for
other baskets. CFMMs use a simple rule for accepting trades: a trade is only valid if the value of a
given function at the post-trade reserves (with a small adjustment to account for fees collected) is
equal to the value at the pre-trade reserves. This function is called the \emph{trading function}
and gives CFMMs their name. A common example of a trading function is the constant product
popularized by Uniswap~\cite{uniswap}, wherein a trade is only accepted if it preserves the product
of the reserve amounts.

CFMMs have quickly become one of the most popular applications of public blockchains, facilitating
several billion dollars of trading volume per day. As DEXs have grown in popularity, so have the
number of CFMMs and assets offered, creating complexity for traders who simply want to maximize
their utility for trading one basket of assets for another. As a result, several ``DEX aggregators"
have emerged to route orders across multiple CFMMs on behalf of users. These aggregators currently
execute several billion dollars per month across all DEXs on Ethereum~\cite{block_aggregators}. At
the same time, CFMM platforms such as Uniswap offer software for routing orders across the subset of
CFMMs they support~\cite{uni_router}.

While the properties of individual CFMMs have been studied
extensively (see, \eg, \cite{angerisImprovedPriceOracles2020, angerisConstantFunctionMarket2021})
it is more
common for users to want to access liquidity on multiple CFMMs to minimize their total trading
costs. In this paper, we study the optimal routing problem for CFMMs. We consider a user who can
trade with multiple CFMMs in order to exchange one basket of assets for another and ask how one
should perform such trades optimally. We show that, in the absence of fixed costs, the optimal
routing problem can be formulated as a convex optimization problem, which is efficiently solvable.
As an (important) sub-case, we demonstrate how to use this problem to identify arbitrage opportunies
on a set of CFMMs. Our framework encompasses the routing problems considered in prior
work~\cite{danos2021global, amm_compositions} and offers solutions in the more general
case where users seek to
trade any basket of assets for any other basket across any set of CFMMs whose trading functions
are concave and not necessarily differentiable.
When including transaction costs, we show that the optimal routing problem is a
mixed-integer convex problem, which permits (potentially computationally intensive) global solutions
as well as approximate solutions using  heuristics based on convex optimization.

\paragraph{Outline.}
We describe the optimal routing problem in~\S\ref{s-orp}, ignoring the
fixed transaction costs.
In~\S\ref{s-opt-cond} we give the optimality conditions for the optimal
routing problem, and give conditions under which the optimal action is to
not trade at all.
We give some examples of the optimal routing problem in~\S\ref{s-examples},
including as a special case the detection of arbitrage in the network.
We give a simple numerical example in~\S\ref{s-num-example}.
In~\S\ref{s-trans-cost} we show how to add fixed transaction costs to the
optimal routing problem, and briefly describe some exact and approximate 
solution approaches.

\section{Optimal routing problem}\label{s-orp}

\paragraph{Network of CFMMs.}
We consider a set of $m$ CFMMs, denoted $i=1, \ldots, m$, each of which
trades multiple tokens from among a universe of $n$ tokens, labeled
$j=1 \ldots, n$.
We let $n_i$ denote the number of tokens that CFMM $i$ trades, with
$2 \leq n_i \leq n$.
We can think of the $m$ CFMMs as the vertices of a network or hypergraph,
with $n$ hyper-edges, each corresponding to one of the $n$ assets,
adjacent to those CFMMs that trade it.
Alternatively we can represent this as a bipartite graph, 
with one group of $m$ vertices the CFMMs, and the other group of $n$
vertices the tokens, with an edge between a CFMM and a token if
the CFMM trades the token. 

\begin{figure}
\centering
\begin{tikzpicture}[thick,
  every node/.style={draw,circle},
  cfmm/.style={fill=black},
  every fit/.style={ellipse,draw,inner sep=-2pt,text width=2cm},
  shorten >= 3pt,shorten <= 3pt
]

\begin{scope}[start chain=going below,node distance=7mm]
\foreach \i in {1,2,...,5}
  \node[cfmm,on chain] (c\i) [label=left: \i] {};
\end{scope}

\begin{scope}[xshift=4cm,yshift=-1cm,start chain=going below,node distance=7mm]
\foreach \i in {1,2,...,3}
  \node[on chain] (t\i) [label=right: \i] {};
\end{scope}

\node [fit=(c1) (c5),label=above:CFMMs] {};
\node [fit=(t1) (t3),label=above:Tokens] {};

\draw (c1) -- (t1);
\draw (c1) -- (t2);
\draw (c1) -- (t3);

\draw (c2) -- (t1);
\draw (c2) -- (t2);

\draw (c3) -- (t2);
\draw (c3) -- (t3);

\draw (c4) -- (t1);
\draw (c4) -- (t3);

\draw (c5) -- (t1);
\draw (c5) -- (t3);

\end{tikzpicture}
\caption{Example CFMM network with $m=5$ CFMMs and $n=3$ tokens.}
\label{f-network}
\end{figure}
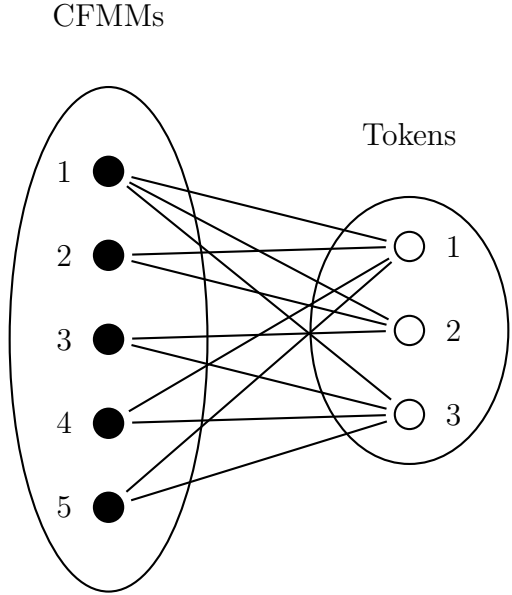

A simple example is illustrated as a bipartite graph in figure~\ref{f-network}.
In this network there are $m=5$ CFMMs, which trade subsets of $n=3$ tokens.
CFMM~1 trades all three tokens, so $n_1=3$; the remaining 4 CFMMs each trade
pairs of tokens, so $n_2= \cdots = n_5=2$.

\paragraph{Global and local token indices.}
We use multiple indices to label the tokens.  The \emph{global} 
index uses the token labels $1, \ldots, n$.  
CFMM $i$, which trades a subset of $n_i$ of the $n$ tokens, 
has its \emph{local} token index, $j=1, \ldots, n_i$.
To link the global and local indexes, we introduce matrices $A_i \in
\reals^{n \times n_i}$, $i=1\ldots, m$, where $(A_i)_{jk} = 1$
if token $k$ in the CFMM $i$'s local index corresponds to global 
token index $j$, while $(A_i)_{jk} =0$ otherwise.

As an example, consider the simple network shown 
in figure~\ref{f-network}. CFMM~1 trades all three tokens, with its local
indexing identical to the global indexing, so $A_1$ is the $3\times 3$ 
identity matrix.  As a more interesting example, consider CFMM~4, 
which trades two tokens, labeled
1 and 2 in the local indexing, but 1 and 3 in the global indexing,
so 
\[
A_4 = \begin{bmatrix}
    1 & 0\\
    0 & 0\\
    0 & 1\\
 \end{bmatrix}.
\]

\paragraph{CFMM tendered and received baskets.}
For CFMM $i$ we denote the tendered and received baskets as
$\Delta_i,~\Lambda_i \in \reals_+^{n_i}$.  These are quantities of
tokens (in the local indexing) that we give to, and receive from,
CFMM $i$ in a proposed trade.   

\paragraph{CFMM trading semantics.}
The proposed trade $(\Delta_i,\Lambda_i)$ is valid and accepted by CFMM
$i$ provided 
\[
\phi_i(R_i+\gamma \Delta_i - \Lambda_i) = \phi_i(R_i),
\]
where $\phi: \reals^{n_i}_+ \to \reals$ is the trading function,
$R_i \in \reals_+^{n_i}$ are the current reserves, and 
$\gamma_i \in (0,1]$ is the trading fee for CFMM $i$.
We will assume that the trading functions
$\phi_i$ are concave and increasing functions.
(For much more detail, see \cite{angerisConstantFunctionMarket2021}.)

Examples of trading functions are the sum function, with
\[
\phi_i(R)= R_1+ \cdots + R_{n_i},
\]
the product or geometric mean
function 
\[
\phi_i(R)= (R_1 \cdots R_{n_i})^{1/n_i},
\]
and a generalization, the weighted geometric mean function
\[
\phi_i(R)= R_1^{w_1} \cdots R_{n_i}^{w_{n_i}},
\]
with weights $w>0$, $\ones^Tw=1$,
where $\ones$ is the vector with all components one.

\paragraph{Network trade vector.}
The CFMM tendered and received baskets $\Delta_i$ and $\Lambda_i$ 
can be mapped to the global indices as the $n$-vectors
$A_i \Delta_i$ and $A_i \Lambda_i$,
which give the numbers of tokens tendered to and received from CFMM $i$,
using the global indices for the tokens.
Summing the difference of received and tendered tokens
over all CFMMs we obtain the \emph{network trade vector}
\[
\Psi = \sum_{i=1}^m A_i (\Lambda_i - \Delta_i).
\]
This is an $n$-vector, which gives the total net number of tokens 
tendered to the CFMMs in the network.
We interpret $\Psi$ as the net network trade vector.  If 
$\Psi_k \geq 0$, it means that we do not need to tender 
any of token $k$ to the network.  This does not mean we 
that do not trade token $k$; it only means that in our 
trading with the CFMMs, we receive at least as much
as token $k$ as we tender.

\paragraph{Network trade utility.}
We introduce a utility function $U:\reals^n \to \reals \cup \{-\infty\}$
that gives the utility of a trade $\Psi$ to a trader as $U(\Psi)$.
We will assume that $U$ is concave and increasing.
Infinite values of $U$ are used to impose constraints;
we consider a proposed trade with $U(\Psi) = -\infty$ as unacceptable.
As an important example, consider the constraint
$\Psi + h \geq 0$, where $h \in \reals_+^n$ is the vector of a user's
current holdings of tokens.  This constraint specifies that the post-trade holdings
$\Psi+ h$ should be nonnegative, \ie, we cannot enter into a trade
that requires more of any token than we currently have on hand.
To express this in the utility function, we define $U(z)$ to be 
$-\infty$ when $z+g \not\geq 0$.  
(This modified utility is also concave and increasing.)

There are many possible choices for $U$.  Perhaps the simplest
is the linear utility function
$U(z) = \pi^Tz$, with $\pi \in \reals_{++}$, where we interpret $\pi_i$ as
the trader's internal or private value or price of token $i$.
Several other choices are discussed
in~\cite[\S5.2]{angerisConstantFunctionMarket2021}.

%

\paragraph{Optimal routing problem.}
We wish to find a set of valid trades that maximizes the trader's utility.
This optimal routing problem can be expressed as
\BEQ\label{e-or-prob}
\begin{aligned}
    & \text{maximize} &&  U(\Psi)\\
    & \text{subject to} && {\textstyle \Psi = \sum_{i=1}^m 
A_i(\Lambda_i - \Delta_i)}\\
    &&& \phi_i(R_i + \gamma_i \Delta_i - \Lambda_i) = \phi_i(R_i), \quad i=1, \dots, m\\
    &&& \Delta_i \geq 0, \quad \Lambda_i \geq 0, \quad i=1, \ldots, m.
\end{aligned}
\EEQ
The variables here are $\Psi \in \reals^n$, 
$\Lambda_i \in \reals_+^{n_i}$,
$\Delta_i \in \reals^{n_i}_+$, $i=1, \ldots, m$.
The data are the utility function $U$, the global-local matrices 
$A_i$, and those associated with the CFMMs, the 
trading functions $\phi_i$,
the trading fees $\gamma_i$, and the reserve 
amounts $R_i \in \reals^{n_i}_+$. 

\paragraph{Convex optimal routing problem.}
Unless the trading functions are affine, the optimal routing 
problem \eqref{e-or-prob} is not a convex optimization problem.
We can, however, form an equivalent problem that is convex.
To do this we replace the equality constraints 
with inequality constraints:
\BEQ\label{e-or-prob-cvx}
\begin{aligned}
    & \text{maximize} &&  U(\Psi)\\
    & \text{subject to} && {\textstyle 
\Psi =\sum_{i=1}^m A_i(\Lambda_i-\Delta_i)}\\
    &&& \phi_i(R_i + \gamma_i \Delta_i - \Lambda_i) \ge \phi_i(R_i), \quad i=1, \dots, m\\
    &&& \Delta_i \geq 0, \quad \Lambda_i \geq 0, \quad i=1, \ldots, m.
\end{aligned}
\EEQ
This problem is evidently convex~\cite[\S5.2.1]{cvxbook} and can be readily solved.

We will show that any solution of 
\eqref{e-or-prob-cvx} is also a solution of \eqref{e-or-prob}.
This the same as showing that for any solution of \eqref{e-or-prob-cvx},
the inequality constraints hold with equality.
Suppose that $\Delta_i^\star$ and $\Lambda_i^\star$ are 
feasible for \eqref{e-or-prob-cvx}, but 
$\phi_k(R_k + \gamma_k \Delta^\star_k - \Lambda^\star_k) > \phi_k(R_k)$
for some $k$.
This means we can find $\tilde \Lambda > \Lambda^\star_k$ which satisfies
$\phi_k(R_k + \gamma_k \Delta^\star_k - \tilde \Lambda) \geq \phi_k(R_k)$.
The associated trade vector $\tilde \Psi = \Psi^\star + A_k \Lambda_k$
satisfies $\tilde \Psi \geq \Psi^\star$, 
$\tilde \Psi \neq \Psi^\star$, \ie, at least one component of $\tilde \Psi$
is larger that the corresponding component of $\Psi^\star$.
It follows that $U(\tilde \Psi) > U(\Psi^\star)$, so $\Psi^\star$ is not optimal
and therefore cannot be a solution.
We conclude that any solution of \eqref{e-or-prob-cvx} satisfies the
inequality constraints as equality, and so is optimal for 
\eqref{e-or-prob}.

A similar statement holds in the case when $U$ is only nondecreasing, and not
(strictly) increasing.  In this case we can say that there is a solution of
\eqref{e-or-prob-cvx} that is optimal for \eqref{e-or-prob}.
One simple method to find such a solution is to solve \eqref{e-or-prob-cvx} 
with objective $U(\Psi)+\epsilon \ones^T\Psi$,
where $\epsilon$ is small and positive.  The objective for this problem is
increasing.  In principle we can let $\epsilon$ go to zero to recover a 
solution of the original problem; in practice choosing a single small value
of $\epsilon$ works.

\paragraph{Consequences of convexity.}
Since the problem \eqref{e-or-prob-cvx} is convex,
it can be reliably and quickly solved, even for large
problem instances~\cite[\S1]{cvxbook}. 
Domain specific languages for convex optimiztion such as
CVXPY~\cite{diamond2016cvxpy, cvxpy_rewriting} 
or JuMP \cite{dunningJuMPModelingLanguage2017}
can be used to specify the optimal routing 
problem in just a few lines of code;
solvers such as ECOS~\cite{ecos}, SCS~\cite{scs}, or Mosek~\cite{mosek} 
can be used to solve the problem.

\paragraph{Implementation.} The solution to~\eqref{e-or-prob-cvx} provides 
the optimal values $(\Delta_i,\Lambda_i)$ that one must trade with 
each CFMM $i$. Note that we do not explicitly consider the ideal execution 
of these trades, as these will depend on the semantics of the underlying 
blockchain that the user is interacting with. For example, users may wish 
to execute trades in a particular sequence, starting with an initial 
portfolio of assets and updating their asset composition after each 
transaction until the final basket is obtained. Users may alternatively use 
flashloans~\cite{flashloans_aave} to atomically perform all trades in a 
single transaction, netting out all trades before repaying the loan at the 
end of the transaction.

\section{Optimality conditions}\label{s-opt-cond}

Assuming that $U$ is differentiable (which it need not be),
the optimality conditions of problem~\eqref{e-or-prob-cvx} are 
feasibility, and the dual conditions
\BEQ\label{e-marg-util} 
\nabla U(\Psi) = \nu,
\EEQ
and
\BEQ\label{e-condi}
 \gamma_i\lambda_i \nabla\phi_i(R_i + 
\gamma_i\Delta_i - \Lambda_i) \le A_i^T\nu \le
\lambda_i\nabla\phi_i(R_i + \gamma_i \Delta_i - \Lambda_i), 
\quad i=1, \ldots, m,
\EEQ
where $\nu \in \reals^n$, $\lambda_i \in \reals_+$,
$i=1,\ldots, m$ are the Lagrange multipliers. 
(We derive these conditions in appendix~\ref{a-opt-condition}.)

The first condition \eqref{e-marg-util} has a very simple interpretation: 
$\nu$ is the vector of marginal utilities of the tokens.
The second set of conditions~\eqref{e-condi}
has a simple interpretation that is very similar to the one given
in~\cite[\S5.1]{angerisConstantFunctionMarket2021} for a 
single CFMM.
The term $\nabla \phi_i(R_i + \gamma_i\Delta_i - \Lambda_i)$, which we 
will write as $P_i \in \reals^{n_i}_+$, can be interpreted as the 
unscaled prices of the tokens that CFMM $i$
trades~\cite[\S2.5]{angerisConstantFunctionMarket2021}, in the local indexing. 
Plugging \eqref{e-marg-util} into~\eqref{e-condi},we get
\BEQ\label{e-opt-condition}
\gamma_i\lambda_iP_i \le A_i^T\nabla U(\Psi) \le \lambda_iP_i, \quad i=1, \dots, m.
\EEQ
The middle term is the vector of marginal utilities of the tokens
that CFMM $i$ trades, in the local indexing.
These marginal utilities must lie between a multiple of the discounted prices,
scaled, and the same prices, adjusted by $\gamma_i$.

We can also recognize the conditions~\eqref{e-condi} as those for 
the problem of finding an optimal trade for CFMM $i$ alone,
with the linear utility $U_i(z) = \pi_i^T z$,
where $\pi_i = A_i^T \nabla U(\Psi)$.
(See~\cite[\S5.2]{angerisConstantFunctionMarket2021}.)
This is very appealing: it states that the optimal trades for
the network of CFMMs are also optimal for each CFMM separately,
when they use a linear utility, with prices equal to the marginal utility of
the overall trade.

\paragraph{Non-differentiable utility.}
When $U$ is not differentiable, the optimality conditions are the same,
but in \eqref{e-marg-util} 
we substitute a supergradient $g \in -\partial (-U)(\Psi)$ for the 
gradient $\nabla U(\Psi)$, where $\partial$ denotes the subdifferential.
When $U$ is not differentiable at $\Psi$, there are multiple such $g$'s,
which we can consider to be multiple marginal utilities.  The optimality condition
is that \eqref{e-opt-condition} hold, with $\nabla U(\Psi)$ replaced with $g$,
for any $g \in -\partial (-U)(\Psi)$.

\paragraph{No-trade condition.}
From the optimality conditions we can derive conditions under which the optimal 
trades are zero, \ie, we should not trade.
We will assume that $U(0) > -\infty$, \ie, $\Psi=0$ is feasible; if this is not the 
case, then evidently $\Psi=0$ is not optimal.
The zero trade 
\[
\Psi=0, \qquad \Delta_i = \Lambda_i = 0, \quad i=1, \ldots, m,
\]
is feasible for \eqref{e-opt-condition}, so the optimality condition is
\BEQ\label{e-no-trade}
\gamma_i\lambda_i P_i \le A_i^T\nabla U(0) \le \lambda_iP_i, 
\quad i=1, \ldots, m,
\EEQ
where $P_i = \nabla \phi_i(R_i)$ is the unscaled price of CFMM $i$ at the current 
reserves $R_i$.
This condition is a generalization of the no-trade condition given 
in~\cite[\S5.1]{angerisConstantFunctionMarket2021} for one CFMM,
to the case where there are multiple CFMMs.
When $U$ is not differentiable, we replace $\nabla U(0)$ with 
a supergradient $g \in -\partial (-U)(0)$.

\section{Examples} \label{s-examples}

\paragraph{Linear utility.}
Consider the linear utility $U(z)=\pi^Tz$, with $\pi >0$.
In this case the optimal routing problem is separable across the CFMMs,
since the objective is a sum of functions of $\Lambda_i-\Delta_i$, and 
constraints are also only on $(\Lambda_i,\Delta_i)$.
It follows that we can solve the optimal routing problem \eqref{e-or-prob-cvx}
by solving $m$ single CFMM problems independently, using linear utilities
with prices given by $\pi$.

\paragraph{Liquidating a basket of tokens.}\label{s-liquidation}
Suppose we start with an initial holding (basket) 
of tokens $h^\text{init} \in \reals_+^n$
and wish to convert them all to token $k$.
We use the utility function
\[
U(z) = \left\{ \begin{array}{ll} z_k & h^\text{init} + z \geq 0\\
-\infty & \mbox{otherwise.} \end{array}\right.
\]
(This utility is nondecreasing but not increasing.)

A special case of this problem is converting one token into another, \ie, when
\[
h^\mathrm{init} = t e_j,
\]
where $t \ge 0$ and $e_j \in \reals^n$ is the $j$th unit vector. In this special
case, we can write the optimal value of the optimization problem, which we will call
$u(t)$, as a function of $t$. It is not hard to show that $u$ is nonnegative
and increasing. This function is also concave as it is the partial maximization
of a concave function (over all variables except $t$). We show an example instance
of this problem, along with an associated function $u$, in~\S\ref{s-examples}.

\paragraph{Purchasing a basket of tokens.} This is the opposite of 
liquidating a basket of tokens.
Here too we start with initial token holdings 
$h^\text{init} \in \reals_+^n$, and end up with the 
holdings $h^\text{init} +  \Psi$.
Let $h^\text{des} \in \reals_+^n$ be our target basket; 
we wish to end up with the largest possible multiple of this 
basket.  Let $\mathcal K \subseteq \{1, \ldots, n\}$
denote the set of indices for which $h^\text{des}_i>0$, \ie,
the indexes associated with tokens in the desired basket.
We seek to maximize the value of $\alpha$ for which
$h^\text{init} + \psi \geq \alpha h^\text{des}$.
To do this we use the utility function
\[
U(z) = \left\{ \begin{array}{ll} \min_{i \in \mathcal K} (h_i^\text{init}+\Psi_i)/
h_i^\text{des} & h^\text{init} + z \geq 0\\
-\infty & \mbox{otherwise.} \end{array}\right.
\]

\paragraph{Arbitrage detection.}\label{s-arbitrage}
An arbitrage is a collection of valid CFMM trades with $\Psi \geq 0$ and 
$\Psi \neq 0$, \ie, a set of trades for the CFMMs in which we tender no tokens,
but receive a positive amount of at least one token.
The optimal routing problem can be used to find an
arbitrage, or certify that no arbitrage exists.

Consider any $U$ that is increasing, 
with domain $\{\Psi \mid U(\Psi)>-\infty \} = \reals_+^n$.
Evidently there is an arbitrage if and only if there is a 
nonzero solution of the routing problem, which is the same as
$U(\Psi^\star) > U(0)$, where $\Psi^\star$ is optimal.
So by solving this optimal routing problem, we can find
an arbitrage, if one exists.

\paragraph{No-arbitrage condition.}
Using the version of \eqref{e-no-trade} for nondifferentiable $U$ we can derive 
conditions under which there is no arbitrage.
We consider the specific utility function 
\[
U(\Psi) = \left\{ \begin{array}{ll} \ones^T \Psi & \Psi \geq 0\\
-\infty & \mbox{otherwise},
\end{array}\right.
\]
where $\ones$ denotes the vector with all entries one.
This utility is the total number of tokens received, when they are nonnegative.
Its supergradient at $0$ is
\[
-\partial (-U)(0) = [1,\infty)^n = \{ g\mid g \geq \ones \}.
\]
The condition \eqref{e-no-trade} becomes:  There exists $g \geq \ones$,
and $\lambda_i \geq 0$, for which 
\[
\gamma_i\lambda_i P_i \le A_i^T g \le \lambda_iP_i, \quad i=1, \ldots, m.
\]
By absorbing a scaling of $g$ into the $\lambda_i$, we can say that 
$g>0$ is enough.

This makes sense: it states that there is no arbitrage if we can 
assign a set of positive prices (given by $g$) to the tokens,
for which no CFMM would trade.
In the (unrealistic) case when $\gamma_i=1$, \ie, there is no trading cost,
the no arbitrage condition is that there exists a global set of prices for the tokens, $g$,
consistent with the local prices of tokens given by $\lambda_iP_i$.

\section{Numerical example}\label{s-num-example}
The Python code for the numerical example we present here is available at
\begin{center}
\texttt{https://github.com/angeris/cfmm-routing-code}.
\end{center}
The optimization problems are formulated and solved using CVXPY \cite{cvxpy_rewriting}.
A listing of the core of the code is given in appendix~\ref{a-code}.

\paragraph{Network.} We consider the network of 5 CFMMs and 3 tokens
shown in figure~\ref{f-network}.  The trading functions, Fee parameters,
and reserves are listed in table~\ref{t-cfmm-example}.

\begin{table}
\centering
\begin{tabular}{l|l|l|l}
    {\bf CFMM} & {\bf Trading function $\phi_i$} &
 {\bf Fee parameter $\gamma_i$} & {\bf Reserves} $R_i$
\\
    \hline
    1 & Geometric mean, $w=(3, 2, 1)$ & 0.98 & (3, .2, 1)\\
    \hline
    2 & Product & 0.99 & (10, 1) \\
    \hline
    3 & Product & 0.96 & (1, 10) \\
    \hline
    4 & Product & 0.97 & (20, 50)  \\
    \hline
    5 & Sum & 0.99 & (10, 10) 
\end{tabular}
\caption{CFMM attributes.}
\label{t-cfmm-example}
\end{table}

\paragraph{Problem and utility.}
We wish to trade an amount $t\geq 0$ of token~1 for the maximum possible
amount of token~3.
This is a special case of the problem of liquidating a basket of tokens,
as described in~\S\ref{s-liquidation}, with initial holdings $h^\text{init}=t e_1$.
The utility function is $U(z) = z_3$, provided $z+h^\text{init} \geq 0$,
and $U(z) = -\infty$ otherwise.
We let $u(t)$ denote the maximum amount of token~3 we can obtain
from the network when we tender token~1 in the amount $t$.

\begin{figure}
    \centering
    \includegraphics[width=0.7\textwidth]{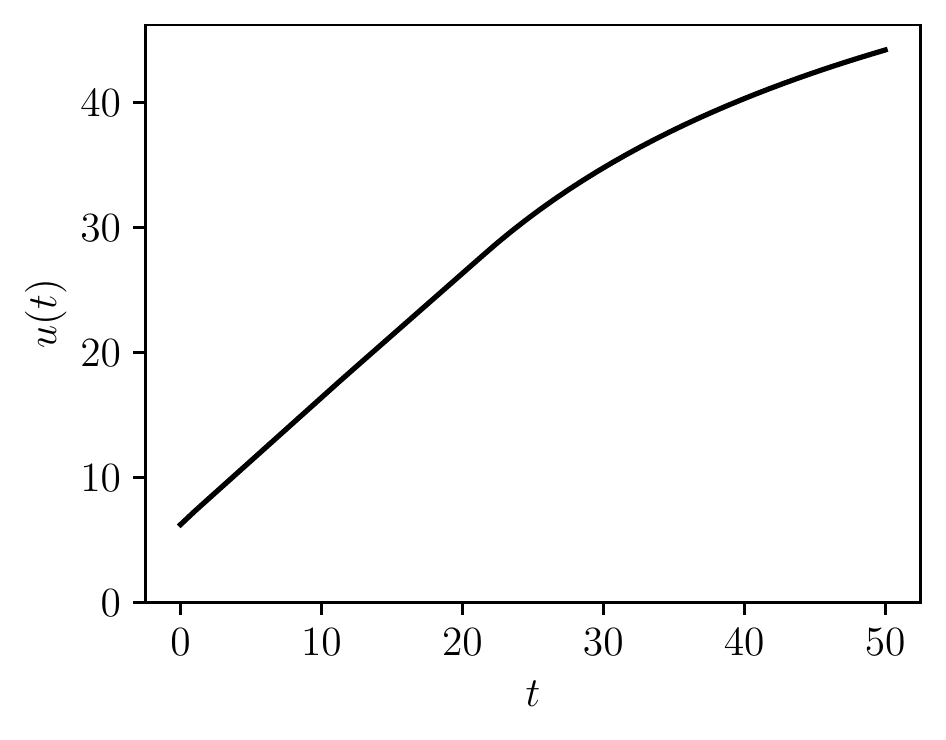}
    \caption{Plot of $u(t)$, the maximum amount of token~3 we obtain
when we tender the amount $t$ of token~1.}
    \label{f-ut}
\end{figure}

\paragraph{Results.}
We solve the optimal routing problem for many values of $t$, from $t=0$ to $t=50$.
The amount of token $3$ we obtain in shown in
figure~\ref{f-ut}. 
We see that $u(0) > 0$, which means there is an arbitrage 
in this network; indeed, there is a set of trades that requires giving zero net 
tokens to the network, but we receive an amount around 7 of token~3.

\begin{figure}
    \centering
    \includegraphics[width=.98\textwidth]{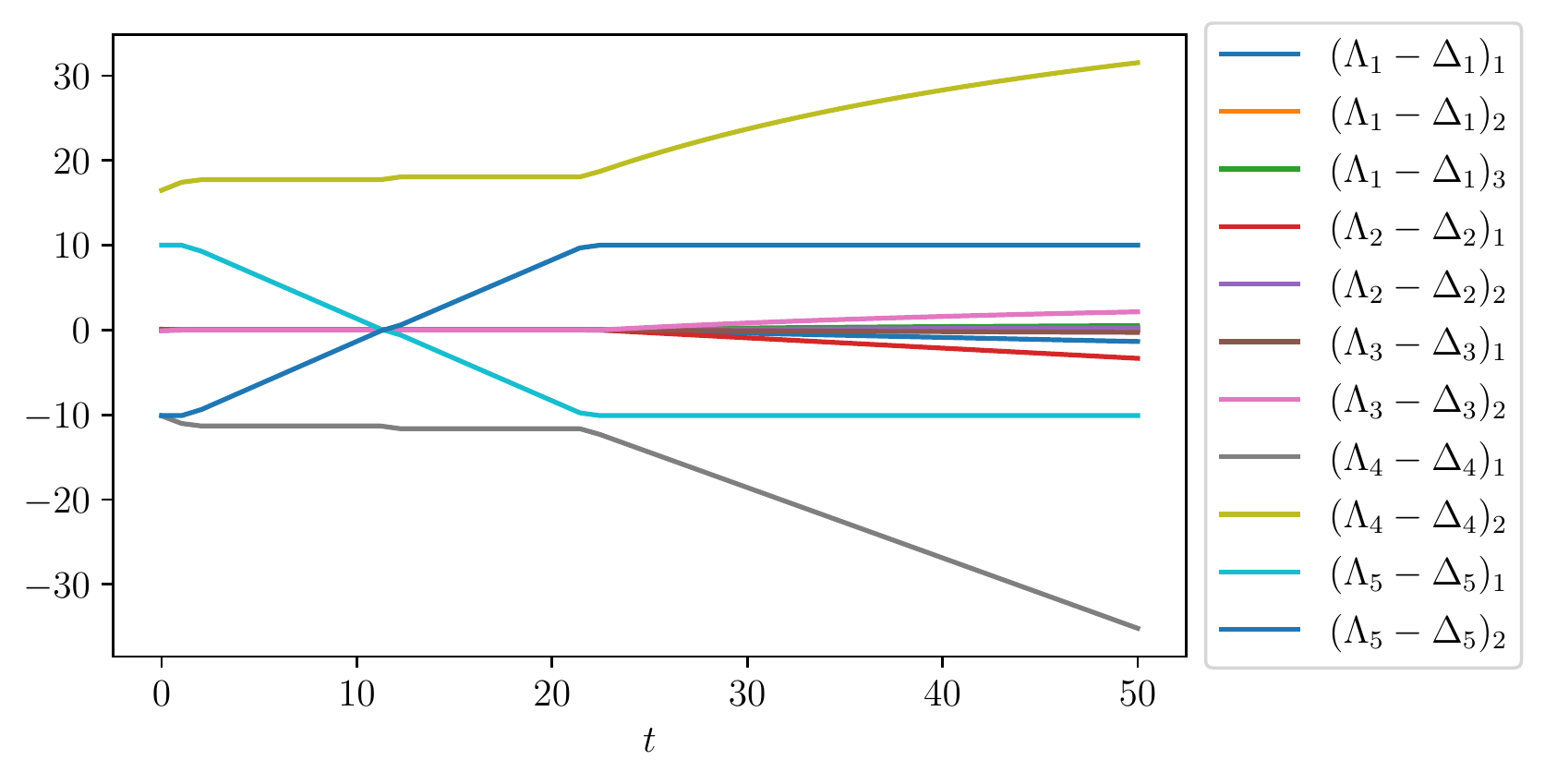}
    \caption{Optimal trades versus $t$, the amount of token~1 tendered.}
\label{f-optimal-trades}
\end{figure}

The associated optimal trades are shown in figure~\ref{f-optimal-trades}.
We can see many interesting phenomena here.  At $t=0$ we see the arbitrage trades,
indeed, the arbitrage trades that yield the largest amount of token~3.
Several asset flows reverse sign as $t$ varies.
For example, for $t<11$, we receive token~1 from CFMM~1, whereas 
for $t>11$, we tender token~1 to CFMM~1.
We can also see that the sparsity pattern of the optimal trades changes with $t$.

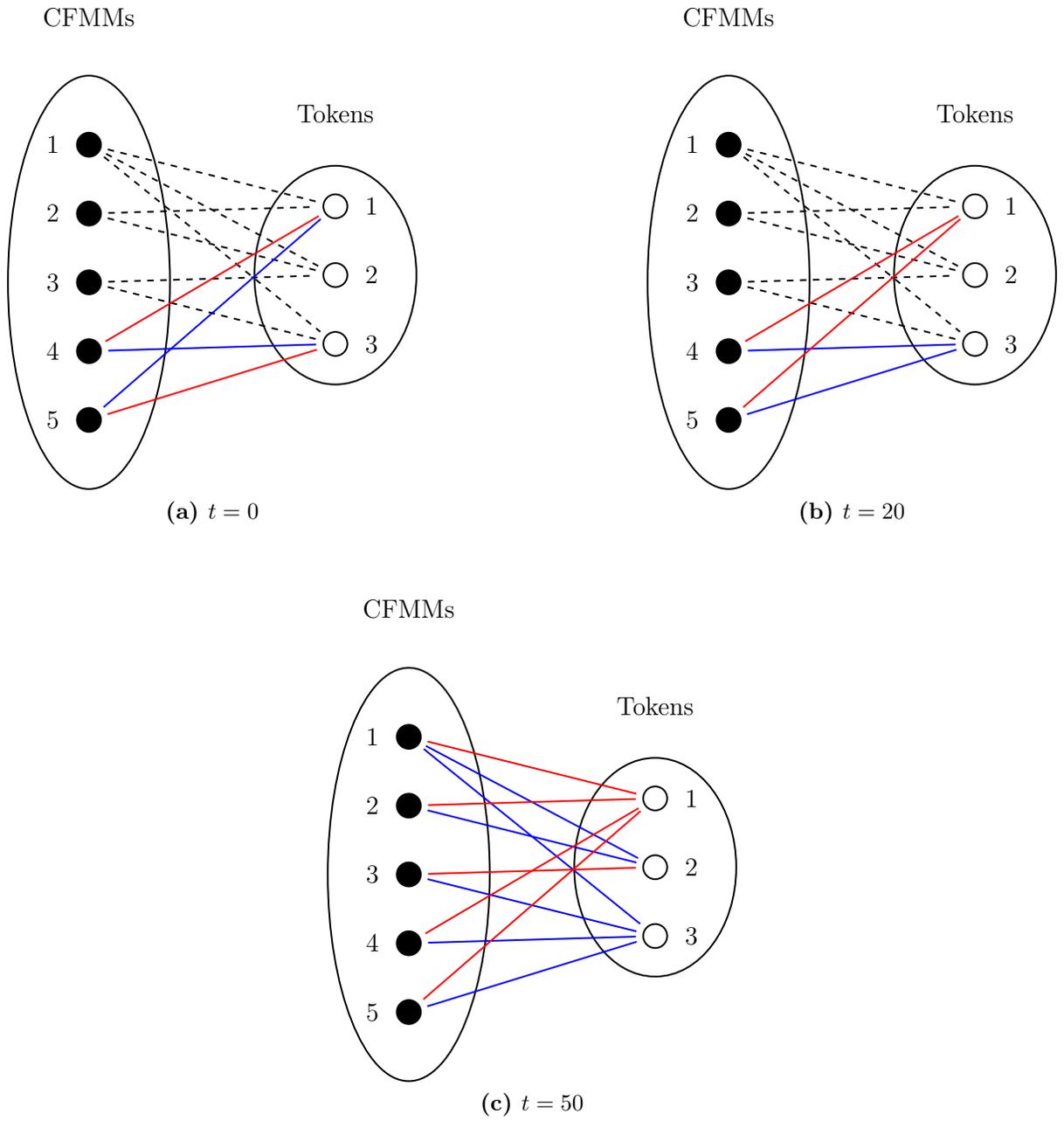
\begin{figure}
\centering
\subfloat[$t=0$]{\label{f-t-.1}
\resizebox{6.5cm}{!}{
\begin{tikzpicture}[thick,
  every node/.style={draw,circle},
  cfmm/.style={fill=black},
  every fit/.style={ellipse,draw,inner sep=-2pt,text width=2cm},
  shorten >= 3pt,shorten <= 3pt
]

\begin{scope}[start chain=going below,node distance=7mm]
\foreach \i in {1,2,...,5}
  \node[cfmm,on chain] (c\i) [label=left: \i] {};
\end{scope}

\begin{scope}[xshift=4cm,yshift=-1cm,start chain=going below,node distance=7mm]
\foreach \i in {1,2,3}
  \node[on chain] (t\i) [label=right: \i] {};
\end{scope}

\node [fit=(c1) (c5),label=above:CFMMs] {};
\node [fit=(t1) (t3),label=above:Tokens] {};

\draw[dashed] (c1) -- (t1);
\draw[dashed] (c1) -- (t2);
\draw[dashed] (c1) -- (t3);

\draw[dashed] (c2) -- (t1);
\draw[dashed] (c2) -- (t2);

\draw[dashed] (c3) -- (t2);
\draw[dashed] (c3) -- (t3);

\draw[red] (c4) -- (t1);
\draw[blue] (c4) -- (t3);

\draw[blue] (c5) -- (t1);
\draw[red] (c5) -- (t3);

\end{tikzpicture}
}
}
\hfill
\subfloat[$t=20$]{
\resizebox{6.5cm}{!}{
\begin{tikzpicture}[thick,
  every node/.style={draw,circle},
  cfmm/.style={fill=black},
  every fit/.style={ellipse,draw,inner sep=-2pt,text width=2cm},
  shorten >= 3pt,shorten <= 3pt
]

\begin{scope}[start chain=going below,node distance=7mm]
\foreach \i in {1,2,...,5}
  \node[cfmm,on chain] (c\i) [label=left: \i] {};
\end{scope}

\begin{scope}[xshift=4cm,yshift=-1cm,start chain=going below,node distance=7mm]
\foreach \i in {1,2,3}
  \node[on chain] (t\i) [label=right: \i] {};
\end{scope}

\node [fit=(c1) (c5),label=above:CFMMs] {};
\node [fit=(t1) (t3),label=above:Tokens] {};

\draw[dashed] (c1) -- (t1);
\draw[dashed] (c1) -- (t2);
\draw[dashed] (c1) -- (t3);

\draw[dashed] (c2) -- (t1);
\draw[dashed] (c2) -- (t2);

\draw[dashed] (c3) -- (t2);
\draw[dashed] (c3) -- (t3);

\draw[red] (c4) -- (t1);
\draw[blue] (c4) -- (t3);

\draw[red] (c5) -- (t1);
\draw[blue] (c5) -- (t3);

\end{tikzpicture}
}
}

\subfloat[$t=50$]{
\resizebox{6.5cm}{!}{
\centering
\begin{tikzpicture}[thick,
  every node/.style={draw,circle},
  cfmm/.style={fill=black},
  every fit/.style={ellipse,draw,inner sep=-2pt,text width=2cm},
  shorten >= 3pt,shorten <= 3pt
]

\begin{scope}[start chain=going below,node distance=7mm]
\foreach \i in {1,2,...,5}
  \node[cfmm,on chain] (c\i) [label=left: \i] {};
\end{scope}

\begin{scope}[xshift=4cm,yshift=-1cm,start chain=going below,node distance=7mm]
\foreach \i in {1,2,3}
  \node[on chain] (t\i) [label=right: \i] {};
\end{scope}

\node [fit=(c1) (c5),label=above:CFMMs] {};
\node [fit=(t1) (t3),label=above:Tokens] {};

\draw[red] (c1) -- (t1);
\draw[blue] (c1) -- (t2);
\draw[blue] (c1) -- (t3);

\draw[red] (c2) -- (t1);
\draw[blue] (c2) -- (t2);

\draw[red] (c3) -- (t2);
\draw[blue] (c3) -- (t3);

\draw[red] (c4) -- (t1);
\draw[blue] (c4) -- (t3);

\draw[red] (c5) -- (t1);
\draw[blue] (c5) -- (t3);

\end{tikzpicture}
}
}
\caption{Optimal trades as $t$ varies.
A red edge means the CFMM is receiving
tokens; a blue edge means the CFMM is tendering tokens.}
\label{f-two-coin}
\end{figure}

We illustrate the changing signs and sparsity of the optimal trades
in figure~\ref{f-two-coin},
which shows the optimal trades for $t=0$, $t=20$, and $t=50$.
We plot whether each token is tendered to or received from each CFMM
using color coded edges.
A red edge connecting a token to a CFMM means that the CFMM is receiving
this token, while a blue edge denotes that the CFMM is tendering
this token. A dashed edge denotes that the CFMM neither tenders nor 
receives this token.

Even in this very simple example, the optimal trades are not obvious
and involve trading with and among multiple CFMMs. For larger networks,
the optimal trades are even less obvious.

\clearpage
\section{Fixed transaction costs}\label{s-trans-cost}
Our optimal routing problem \eqref{e-or-prob-cvx}
includes the trading costs built into CFMMs, via the parameters $\gamma_i$.
But it does not include the small fixed cost associated with any trade.
In this section we explore how these fixed transaction costs can be incorporated
into the optimal routing problem.

We let $q_i \in \reals_+$ denote the fixed cost of executing a trade with CFMM $i$,
denominated in some numeraire.
We pay this whenever we trade, \ie, $\Lambda_i - \Delta_i \neq 0$.
We introduce a new set of Boolean
variables into the problem, $\eta \in \{0, 1\}^m$, with $\eta_i = 1$ if a
nonzero trade is made with CFMM $i$, and $\eta_i = 0$ otherwise, so
the total fixed transaction cost is $q^T\eta$.
We assume that there is a known maximum size of a
tendered basket with CFMM $i$, which we will denote $\dmax_i \in \reals_+^{n_i}$.
We can then express the problem of maximizing the utility minus the
fixed transaction cost as
\begin{equation}\label{eq:micp}
\begin{aligned}
    & \text{maximize} &&  U(\Psi) - q^T\eta\\
    & \text{subject to} && {\textstyle \Psi =\sum_{i=1}^m A_i(\Lambda_i - \Delta_i)}\\
    &&& \phi_i(R_i + \gamma_i \Delta_i - \Lambda_i) \ge \phi_i(R_i), \quad i=1, \dots, m\\
    &&& 0 \le \Delta_i \le \eta_i\dmax_i, \quad \Lambda_i \ge 0, \quad i=1, \dots, m,\\
    &&& \eta \in\{0, 1\}^m,
\end{aligned}
\end{equation}
where the variables are $\Psi$, $\Lambda_i$, $\Delta_i$, $i=1, \ldots, m$,
and $\eta \in \{0,1\}^m$.

The optimal routing problem with fixed costs \eqref{eq:micp}
is a mixed-integer convex program (MICP).
It can be solved exactly, possibly with great computational effort,
using global optimization methods, with MICP solvers such as
Mosek~\cite{mosek} or Gurobi~\cite{gurobi}. 
When $m$ is small (say, under ten or so), it can be practical to 
solve it by brute force, by solving the convex problem we obtain 
for each of the $2^m$ feasible values of $\eta$.

\paragraph{Approximate solution methods.} 
Many approximate methods have the speed of convex optimization, 
and (often) produce good approximate solutions.
For example we can solve the relaxation of \eqref{eq:micp} obtained
by replacing the constraints on $\eta$ to $\eta \in[0,1]^m$ (which 
gives a convex optimization problem).
After that we set a threshold $t\in (0,1)$ for the relaxed optimal 
values $\eta^\text{rel}$,
and take $\eta_i = 1$ when $\eta^\text{rel}\geq t$ and
and $\eta_i = 0$ when $\eta^\text{rel}<t$. We fix these values of $\eta$ and
then solve the resulting convex problem.
This could be done for a modest number of values of $t$; we take the solution
found with the largest objective (including the fixed costs $q^T\eta$).

An alternative is a simple randomized method.
We interpret $\eta^\text{rel}_i$ as probabilities, and generate $\eta \in
\{0,1\}$ randomly using these probabilities. We then solve the convex problem
associated with this choice of $\eta$. We can repeat this procedure a modest number of times,
and pick the feasible point with the highest payoff.

\clearpage
\bibliographystyle{alpha}
\bibliography{citations.bib}

\appendix

\clearpage
\section{Derivation of optimality conditions}\label{a-opt-condition}
To derive condition~\eqref{e-condi}, we introduce the Lagrangian of 
the convex optimal routing problem~\eqref{e-or-prob-cvx},
\begin{multline*}
\mathcal{L}(\Psi, \Delta, \Lambda, \nu, \lambda, \eta, \xi) =
-U(\Psi) + \nu^T\left(\Psi - \sum_{i=1}^m A_i(\Lambda_i - \Delta_i)\right) \\
+ \sum_{i=1}^m \lambda_i(\phi_i(R_i) - 
\phi_i(R_i+\gamma_i\Delta_i - \Lambda_i)) - \sum_{i=1}^m 
(\eta_i^T\Delta_i + \xi_i^T\Lambda_i).
\end{multline*}
Here, $\nu \in \reals^n$, $\lambda \in \reals^m_+$, and 
$\eta_i, \xi_i \in \reals^{n_i}_+$, $i=1, \dots, m$. These
are the dual variables or Lagrange multipliers for the equality constraint,
the relaxed inequality constraints, and nonnegativity constraints
for the tendered and received baskets of assets for each CFMM, respectively.

The optimality conditions are (primal) feasibility, and (the dual conditions)
\[
\nabla_\Psi \mathcal{L} = 0, \quad \nabla_\Delta\mathcal{L}=0, \quad 
\nabla_\Lambda\mathcal{L} = 0.
\]
The first of these is
\[
-\nabla U(\Psi) + \nu = 0,
\]
while the latter two are
\[
-A_i^T\nu + \lambda_i\nabla\phi_i(R_i+\gamma_i\Delta_i - \Lambda_i) - \xi_i = 0, 
\qquad A_i^T\nu - \gamma_i\lambda_i\nabla\phi_i
(R_i+\gamma_i\Delta_i - \Lambda_i) - \eta_i = 0.
\]
These can be written as
\[
-A_i^T\nu + \lambda_i\nabla\phi_i(R_i+\gamma_i\Delta_i - \Lambda_i) \ge 0, \qquad A_i^T\nu - \gamma_i\lambda_i\nabla\phi_i(R_i+\gamma_i\Delta_i - \Lambda_i) \ge 0,
\]
which simplifies to~\eqref{e-condi}.

\clearpage
\section{Example code}\label{a-code}
\begin{lstlisting}[language=PythonPlus, style=colorEX]
import numpy as np
import cxvpy as cp

# We assume the following is already defined:
# - `A`, list of matrices A_i
# - `reserves`, list of vectors R_i
# - `fees`, list of fees for each CFMM i
# - `num_tokens`, list of n_i
# - `t`, amount of asset 1 provided to the network

current_assets = np.array([t, 0, 0])

# Build variables
deltas = [cp.Variable(n_i, nonneg=True) for n_i in num_tokens]
lambdas = [cp.Variable(n_i, nonneg=True) for n_i in num_tokens]
psi = cp.sum(A_i @ (L - D) for A_i, D, L in zip(A, deltas, lambdas))

# Objective is to trade t of asset 1 for a maximum amount of asset 3
obj = cp.Maximize(psi[2])

# Reserves after trade
new_reserves = [R + gamma_i*D - L \
    for R, gamma_i, D, L in zip(reserves, fees, deltas, lambdas)]

# Trading function constraints
cons = [
    # Balancer pool with weights 3, 2, 1
    cp.geo_mean(new_reserves[0], p=np.array([3, 2, 1])) \
        >= cp.geo_mean(reserves[0], p=np.array([3, 2, 1])),

    # Uniswap v2 pools
    cp.geo_mean(new_reserves[1]) >= cp.geo_mean(reserves[1]),
    cp.geo_mean(new_reserves[2]) >= cp.geo_mean(reserves[2]),
    cp.geo_mean(new_reserves[3]) >= cp.geo_mean(reserves[3]),

    # Constant sum pool
    cp.sum(new_reserves[4]) >= cp.sum(reserves[4]),
    new_reserves[4] >= 0,

    # Allow all assets at hand to be traded
    psi + current_assets >= 0
]

# Set up and solve problem
prob = cp.Problem(obj, cons)
prob.solve()

print(f"amount of asset 3 received: {psi[2].value}")
\end{lstlisting}
\end{document}